\documentclass[12pt,a4paper,reqno,dvips]{article}%
\usepackage{amsfonts}
\usepackage{bm}
\usepackage[all]{xy}
\usepackage{amsmath}
\usepackage{amssymb}
\usepackage{graphicx}
\SelectTips{cm}{}
\xyoption{dvips}
\newtheorem{theorem}{Theorem}
\newtheorem{corollary}[theorem]{Corollary}
\newtheorem{definition}[theorem]{Definition}

\newtheorem{proposition}[theorem]{Proposition}

\newenvironment{proof}[1][Proof]{\noindent\textbf{#1.} }{\ \rule{0.5em}{0.5em}}
\title{\bf Iterated Differential Forms I: Tensors}
\author{\sc{A.~M.~Vinogradov}\thanks{{\bf e}-{\it mail}: \texttt{vinograd@unisa.it}} and \sc{L.~Vitagliano}\thanks{{\bf e}-{\it mail}: \texttt{luca\_vitagliano@fastwebnet.it}}\\
\small{DMI, Universit\`a degli Studi di Salerno}\\ \small{and INFN, Gruppo collegato di Salerno,}\\
\small{Via Ponte don Melillo, 84084 Fisciano (SA), Italy}}
\begin{document}
\maketitle
\begin{abstract}
This note is the first in a series of short communications
dedicated to general theory and some applications of iterated
differential forms. Both are developed either in the
\textquotedblleft classical\textquotedblright context, or in the
\textquotedblleft quantistic\textquotedblright one, i.e., of
Secondary Calculus (see \cite{v98}, \cite{v01}). Detailed
expositions containing proofs of the announced results will be
appearing in due course.\\ With iterated forms we solve the
problem of a \emph{conceptual foundation of tensor calculus}. In
particular, we show that covariant tensors are differential forms
over a certain graded commutative algebra called the algebra of
iterated differential forms. From one side, this interpretation
extends noteworthy frames of the traditional tensor calculus and
enriches it by numerous new natural operators. On the other side,
it allows various generalizations of tensor calculus, the most
important of which is that to secondary (\textquotedblleft
quantized\textquotedblright) calculus. In particular, this leads
to an unified solution of the secondarization (\textquotedblleft
quantization\textquotedblright) problem (see \cite{v98},
\cite{v01}) for arbitrary tensors.\\ In this communication
the algebra of iterated (differential) forms over an arbitrary
(graded) commutative algebra is defined. It is also shown how
tensors on a (smooth) manifold $M$ are naturally interpreted as
iterated differential forms over the algebra $C^{\infty}(M)$.
\end{abstract}
\newpage According to the original definition a tensor on a manifold $M$ is
a system of suitably indexed functions, called components,
associated with a local chart, which change according to a certain
rule when passing to another chart. The modern definition presents
a tensor to be a $C^{\infty}(M)$--multilinear and
$C^{\infty}(M)$--valued function in a number of variables that are
either vector fields, or differential $1$--forms on $M$. These
definitions are manifestly descriptive. For instance, on the basis
of any of them it is not possible to understand why the natural
exterior differential $d$ is defined on skew--symmetric covariant
tensors, i.e., differential forms, but not on symmetric ones. Or,
why a natural connection, namely, that of Levi--Civita, is
associated with (non--degenerate) symmetric covariant
$2$--tensors, but not with skew--symmetric ones, etc. According to
\cite{v72} the meaning of \textquotedblleft
controvariant\textquotedblright objects of differential calculus,
say, vector fields, or general differential operators, is given by
\emph{functors of differential calculus}. For instance, the
derivation functor $\mathrm{D}$ corresponds to vector fields and
$\mathrm{Diff}_{k}$ to $k$--th order differential operators. On
the other hand, covariant objects are related to representing such
functors objects. For instance, differential $i$--forms and
$k$--jets are elements of modules that represent functors
$\mathrm{D}_{i}$ and $\mathrm{Diff}_{k}$, respectively. Concerning
conceptual definition of a species of covariant tensors the
problem is to attribute them to the module representing a certain
functor of differential calculus. What makes this problem not very
banal is that such a \emph{direct} attribution is not, generally,
possible. This is, for instance, the case of symmetric tensors. On
the contrary, skew-symmetric tensors, i.e., differential forms,
allow such one and the corresponding functors are
$\mathrm{D}_{i}$'s. We overcome this difficulty by looking for the
necessary \emph{direct attribution} not over the original ground
commutative algebra, but over another one, naturally associated
with the former.

Below this idea is realized  by passing to the filtered graded
commutative algebra $A=\Lambda_{0} \subset\Lambda_{1}
\subset\cdots\subset\Lambda_{k}
\subset\cdots\subset\Lambda_{\infty}$, where $A$ is the ground
algebra and $\Lambda_{k}$ is the algebra of differential forms
over the algebra $\Lambda_{k-1}$. This way we respond the
question: what properly  are covariant tensors?

In the subsequent note the first application of the theory
sketched here to Riemannian geometry  will be given. In
particular, the nature of the Levi--Civita connection will be
clarified and answered the question posed above.

\section{Differential calculus over graded algebras}

A pair $\mathcal{G}=(G,\mu),\;G$ being an abelian group and $\mu:G\times
G\longrightarrow\mathbb{Z}_{2}$ a symmetric $\mathbb{Z}$--bilinear map, is a
\emph{grading group}. Below $\Bbbk$ stands for a field of zero--characteristic
and we use $g\cdot h$ for $\mu(g,h)$. The category of unitary, $\mathcal{G}%
$--graded, associative, graded--commutative $\Bbbk$--algebras is denoted by
$\mathbf{Alg}_{\Bbbk}^{\mathcal{G}}$. If $A\in\mathbf{Alg}_{\Bbbk
}^{\mathcal{G}}$ and $a\in A$ is homogeneous, then $|a|\in G$ denotes the
degree of $a$. Recall that graded--commutativity means that $ab=(-1)^{|a|\cdot
|b|}ba$ for all homogeneous elements $a,b\in A$. In the sequel we adopt the
following convention: in the exponent of $(-1)$ the symbol, say, $\sigma$,
denoting a homogeneous element is used as a substitute of $|\sigma|$. For
instance, $(-1)^{a\cdot b}$ means $(-1)^{|a|\cdot|b|}$. The differential
calculus over an algebra $A\in\mathbf{Alg}_{\Bbbk}^{\mathcal{G}}$ is
introduced along the lines of \cite{v72,v81,v89,v96} where the reader will
find further details.

Denote by $\mathbf{Mod}_{A}^{\mathcal{G}}$ the category of $\mathcal{G}%
$--graded $A$--modules and by $\mathrm{D}_{A}:\mathbf{Mod}_{A}^{\mathcal{G}%
}\longrightarrow\mathbf{Mod}_{A}^{\mathcal{G}}$ the functor associating with
$P\in\mathbf{Mod}_{A}^{\mathcal{G}}$ the graded $A$--module $\mathrm{D}%
_{A}(P)$ of $P$\emph{--valued} \emph{derivations} of $A$. $\mathrm{D}_{A}(A)$
is a graded $\Bbbk$--Lie algebra with respect to the graded commutator.
Functors $\mathrm{D}_{k}:\mathbf{Mod}_{A}^{\mathcal{G}}\longrightarrow
\mathbf{Mod}_{A}^{\mathcal{G}}$, $k\in\mathbb{N}$ are defined in
\cite{v72,v81,v89}. Recall that an element $\nabla\in\mathrm{D}_{k}(P)$ may be
viewed as a graded skew--symmetric $P$--valued multi--derivation of $A$ of
multiplicity $k$. In particular, $\mathrm{D}_{1}=\mathrm{D}_{A}$. Denote by
$\Lambda^{k}(A)$ the graded $A$--module of differential $k$--forms over $A$
and by $\Lambda^{k}(A)^{g},\;g\in G$, its homogeneous component of grade $g$.

The direct sum $\Lambda(A)=\bigoplus_{k=0}^{\infty}\Lambda^{k}(A)$ has a
natural structure of an unitary, $\mathcal{G}\oplus\mathbb{Z}$--graded,
associative, graded--commutative $\Bbbk$--algebra which is called the algebra
of \emph{differential forms over} $A$. It is naturally isomorphic to the
$\mathcal{G}$--exterior algebra $\bigwedge^{\ast}\Lambda^{1}(A)$. Thus, any
$\sigma\in$ $\Lambda(A)$ is of the form
\begin{equation}
\sigma=\sum a_{\alpha_{1}\cdots\alpha_{k}}db_{\alpha_{1}}\wedge\cdots\wedge
db_{\alpha_{k}}, \label{Eq0}%
\end{equation}
for some $a_{\alpha_{1}\ldots\alpha_{k},}b_{\alpha_{1}},\ldots b_{\alpha_{k}%
}\in A$.

For any $\sigma\in\Lambda^{k}(A)^{g}$ put $|\sigma|=(g,k)\in
G\oplus\mathbb{Z}$. Then, $\sigma\wedge\rho=(-1)^{\sigma\cdot\rho}\rho
\wedge\sigma$ for any homogeneous elements $\sigma,\rho\in\Lambda(A)$, where
$(g,k)\cdot(h,l)\overset{\mathrm{def}}{=} g\cdot h+[kl]_{2}\in\mathbb{Z}_{2}$.

The \emph{exterior differential} in $\Lambda(A)$ will be denoted by
$d:\Lambda(A)\longrightarrow\Lambda(A)$. It is a graded $\Lambda
(A)$--derivation of bi--degree $|d|=(0,1)\in G\oplus\mathbb{Z}$ and $d^{2}=0$.
So, $(\Lambda(A),d)$ is a $\mathcal{G}\oplus\mathbb{Z}$--graded differential algebra.

The correspondence $A\longmapsto(\Lambda(A),d)$ is a functor from the category
$\mathbf{Alg}_{\Bbbk}^{\mathcal{G}}$ to the category $\mathbf{dAlg}_{\Bbbk
}^{\mathcal{G}\oplus\mathbb{Z}}$ of unitary, $\mathcal{G}\oplus\mathbb{Z}%
$--graded, associative, graded--commutative, differential $\Bbbk$--algebras.
In particular, if $A,A^{\prime}\in\mathbf{Alg}_{\Bbbk}^{\mathcal{G}}$, then a
morphism $\phi:A\longrightarrow A^{\prime}$ is extended to a morphism
$\Lambda(\phi):\Lambda(A)\longrightarrow\Lambda(A^{\prime})$ compatible with
the exterior differential, i.e., $d\circ\Lambda(\phi)=\Lambda(\phi)\circ d$.
If $\sigma\in\Lambda(A)$ is of the form (\ref{Eq0}), then $\Lambda
(\phi)(\sigma)=\sum\phi(a_{\alpha_{1}\cdots\alpha_{k}})d(\phi(b_{\alpha_{1}%
}))\wedge\cdots\wedge d(\phi(b_{\alpha_{k}}))$.

The insertion of a derivation $X \in\mathrm{D}_{A}(A)$ operator will be
denoted by $i_{X}:\Lambda(A)\longrightarrow\Lambda(A)$. It is a $\Lambda
(A)$--derivation of bi--degree $(| X|,-1)$. The Lie derivative along $X
\in\mathrm{D}_{A}(A)$ is defined by $\mathcal{L}_{X}=\lbrack i_{X},d]$. It is
a graded $\Lambda(A)$--derivation of bi--degree $(|X|,0)$ and extends $X$ to
$\Lambda(A)$.

For any $X,Y\in\mathrm{D}_{A}(A)$ the following (graded) commutation relations
hold:
\begin{equation}
\lbrack i_{X},i_{Y}]=[\mathcal{L}_{X},d]=0,\quad\lbrack i_{X},\mathcal{L}%
_{Y}]=i_{[X,Y]},\quad\lbrack\mathcal{L}_{X},\mathcal{L}_{Y}]=\mathcal{L}%
_{[X,Y]}. \label{Eq1}%
\end{equation}

\section{Iterated differential forms}

In what follows we put $\Lambda\equiv\Lambda(A)$. So, $\Lambda$ is an unitary,
graded, associative, graded--commutative algebra and all constructions of the
previous section are applied to $\Lambda$ and so on. This leads to the
following \cite{vv06}

\begin{definition}
Given $k\in\mathbb{\mathbb{N}}$ the \emph{algebra} $\Lambda_{k}$ of
$k$\emph{--times iterated differential forms over }$A$ is defined inductively
as $\Lambda_{k}=\Lambda(\Lambda_{k-1})$, by starting with $\Lambda_{0}=A$. The
exterior differential $d_{k}$ in $\Lambda_{k}$ ($d_{1}=d$) is called
$k$\emph{--th iterated exterior differential}.
\end{definition}

So, $(\Lambda_{k},d_{k})$ is a differential, $\mathcal{G}\oplus\mathbb{Z}%
^{k}$--graded commutative algebra. Natural inclusions $\Lambda_{k-1}%
\subset\Lambda(\Lambda_{k-1})=\Lambda_{k}$ define the filtered algebra
$\Lambda_{\infty}:\Lambda_{0}\equiv A\subset\Lambda_{1}\equiv\Lambda
\subset\Lambda_{2}\subset\cdots\subset\Lambda_{k}\subset\cdots\subset
\Lambda_{\infty}$, where $\Lambda_{\infty}\equiv\bigcup\nolimits_{k}%
\Lambda_{k}$. $\Lambda_{\infty}$ is a $\mathcal{G}\oplus\mathbb{Z}^{\infty}%
$--graded commutative algebra. If $(g,K)\in G\oplus\mathbb{Z}^{k}$, denote by
$\Lambda_{k}^{(g,K)}\subset\Lambda_{k}$ the homogeneous component of grade
$(g,K)$. Put also $\Lambda_{k}^{K} \overset{\mathrm{def}}{=}\bigoplus_{g\in
G}\Lambda_{k}^{(g,K)}\subset\Lambda_{k}$.

\begin{definition}
$\Lambda_{\infty}$ is called the \emph{algebra of iterated forms over} $A$.
\end{definition}

The operator of insertion of $\nabla\in\mathrm{D}_{\Lambda_{k}}(\Lambda_{k})$
into forms $\Lambda(\Lambda_{k})=\Lambda_{k+1}$ over $\Lambda_{k}$ is denoted
by $i_{\nabla}^{(k+1)}\in\mathrm{D}_{\Lambda_{k+1}}(\Lambda_{k+1})$. The Lie
derivative $\mathcal{L}_{\nabla}^{(k+1)}=[i_{\nabla}^{(k+1)},d_{k+1}]
\in\mathrm{D}_{\Lambda_{k+1}}(\Lambda_{k+1})$ extends $\nabla$ to
$\Lambda_{k+1}$. In its turn, $\mathcal{L}_{\nabla}^{(k+1)}$ can be extended
to $\Lambda_{k+2}$ and so on up to $\Lambda_{\infty}$. Hence any $\nabla
\in\mathrm{D}_{\Lambda_{k}}(\Lambda_{k})$ extends to a derivation of
$\Lambda_{\infty}$. This extension will be denoted by the same symbol $\nabla
$. Note that this notation is consistent with the third relation in (\ref{Eq1}).

In particular, the $k$--th differential $d_{k}$ can be extended to a
derivation of $\Lambda_{\infty}$. Moreover, It follows from commutation
relations (\ref{Eq1}) that $d_{j}^{2}=0$ and $[d_{i},d_{j}]=0$ for all $i,j$.
This way $(\Lambda_{\infty},d_{1},\ldots,d_{k},\ldots)$ becomes a multiple
complex. Note also that for a given $k$ the correspondence $A\longmapsto
(\Lambda_{\infty},d_{k})$ is a functor from the category $\mathbf{Alg}_{\Bbbk
}^{\mathcal{G}}$ to the category $\mathbf{dAlg}_{\Bbbk}^{\mathcal{G}%
\oplus\mathbb{Z}^{\infty}}$.

Note that $\Lambda_{k+l}=\Lambda_{l}(\Lambda_{k})$ for any $k,l\in\mathbb{N}$.
Therefore, $\Lambda_{l}(\Lambda_{\infty})=\Lambda_{\infty}$, $l \in\mathbb{N}%
$. Indeed,
\[
\Lambda_{l}(\Lambda_{\infty})=\Lambda_{l}({\textstyle\bigcup\nolimits_{k}%
}\Lambda_{k})={\textstyle\bigcup\nolimits_{k}}\Lambda_{l}(\Lambda
_{k})={\textstyle\bigcup\nolimits_{k}}\Lambda_{k+l}=\Lambda_{\infty}.
\]
We express this fact by saying accordingly that $\Lambda_{\infty}$ \emph{is
}$\Lambda$--\emph{closed}.

\begin{proposition}
\label{PropKappa}There is an isomorphism $\kappa_{(12)}$ of the double
complexes $(\Lambda_{2},d_{1},d_{2})$ and $(\Lambda_{2},d_{2},d_{1})$.
\end{proposition}

\begin{proof}
The proposition is proved via the following commutative diagram
\[
\xymatrix@=9ex{
A \ar[r]^{d_1}  \ar@{=}[d]& \Lambda^1_1\  \ar[d]^{i_{d_2}} \ar@{^{(}->}[r]& \Lambda_1 \ar[r]^{d_2} \ar[d]^{\kappa^\prime}  \ar@/_1.8pc/[dd]_<(0.21){\bigwedge^\ast (i_{d_2})}|\hole& \Lambda_2^1\  \ar@{^{(}->}[r] \ar[d]^{\kappa^{\prime\prime}}  \ar@/_1.8pc/[dd]_<(0.21){\Lambda^1(\kappa^{\prime})}|\hole & \Lambda_2 \ar@<0.5ex>[d]^{\kappa_{(12)}}  \ar@/_1.8pc/[dd]_<(0.21){\bigwedge^\ast \kappa^{\prime\prime}}|\hole \\
A \ar[r]^{d_2} & \Lambda_2^{(0,1)}\ \ar@{^{(}->}[r] \ar@{^{(}->}[dr] & \Lambda_2^{(0,\ast)}\ \ar[r]^{d_1} \ar[dr]_{d^{\prime\prime}} & \Lambda_2^{(1,\ast)}\ \ar@{^{(}->}[r] \ar@{^{(}->}[dr] & \Lambda_2 \ar@<0.5ex>[u] \\
  &                    & \bigwedge_A^\ast (\Lambda_2^{(0,1)})  \ar[u]_{\bigwedge^\prime}& \Lambda^1(\Lambda_2^{(0,\ast)}) \ar[u]_{i_{d^{\prime\prime}}} & \bigwedge_{\Lambda_2^{(0,\ast)}}^\ast (\Lambda_2^{(1,\ast)})   \ar[u]_{\bigwedge^{\prime\prime\prime}}
}
\]
where the following notations have been adopted. $\Lambda_{2}^{(m,\ast
)}=\bigoplus_{l}\Lambda_{2}^{(m,l)}$, $m=0,1$, $\wedge^{\prime}$ is the
product in $\Lambda_{2}$ of elements belonging to $\Lambda_{2}^{(0,1)}%
\subset\Lambda_{2}$ and $\kappa^{\prime}=\wedge^{\prime}\circ\bigwedge^{\ast
}(i_{d_{2}})$. Note that $\Lambda_{2}^{(0,\ast)}$ is a $\mathcal{G}%
\oplus\mathbb{Z}$--graded, unitary, graded--commutative algebra and
$d^{\prime\prime}:\Lambda_{2}^{(0,\ast)}\longrightarrow\Lambda^{1}(\Lambda
_{2}^{(0,\ast)})$ is its first de Rham differential. Moreover $\kappa
^{\prime\prime}=i_{d^{\prime\prime}}\circ\Lambda^{1}(\kappa^{\prime})$.
Finally $\wedge^{\prime\prime\prime}$ is the product in $\Lambda_{2}$ of
elements belonging to $\Lambda_{2}^{(1,\ast)}\subset\Lambda_{2}$ and
$\kappa_{(12)}=\wedge^{\prime\prime\prime}\circ\bigwedge^{\ast}\kappa
^{\prime\prime}$. It is straightforward to see that
\begin{align*}
&  \kappa_{(12)}(d_{1}g_{1}\wedge\cdots\wedge d_{1}g_{p}\wedge d_{2}%
h_{1}\wedge\cdots\wedge d_{2}h_{q}\wedge d_{1}d_{2}\ell_{1}\wedge\cdots\wedge
d_{1}d_{2}\ell_{r})\\
&  =d_{2}g_{1}\wedge\cdots\wedge d_{2}g_{p}\wedge d_{1}h_{1}\wedge\cdots\wedge
d_{1}h_{q}\wedge d_{1}d_{2}\ell_{1}\wedge\cdots\wedge d_{1}d_{2}\ell_{r}%
\end{align*}
for any $g_{1},\ldots,g_{p},h_{1},\ldots,h_{q},\ell_{1},\ldots,\ell_{r}\in A$.
Therefore, $\kappa_{(12)}$ is an involution of $\Lambda_{2}$.\end{proof}

Let $S_{k}$ be the group of permutations of $\{1,\ldots,k\}$.

\begin{corollary}
For any $k\in\mathbb{N}$ and $\sigma\in S_{k}$ there is an isomorphism
$\kappa_{\sigma}$ of the multiple complexes $(\Lambda_{k},d_{1},\ldots,d_{k})$
and $(\Lambda_{k},d_{\sigma(1)},\ldots,d_{\sigma(k)})$.
\end{corollary}

\begin{proof}
The proof is by induction on $k$. Obviously, it is sufficient to prove the
assertion for a transposition $\sigma$. The base of induction, $k=2$, is
provided by the above proposition. Now, suppose the corollary be true for
$\Lambda_{k-1}$. In particular, for any $i,j<k$ there is an isomorphism
$\kappa_{i,j}$ of the complexes $(\Lambda_{k-1},d_{i})$ and $(\Lambda
_{k-1},d_{j})$ which commutes with all the other differentials. By extending
the $d_{i}$'s, $i<k$, as derivations to $\Lambda_{k}$ one finds that the
complexes $(\Lambda_{k},d_{i})$ and $(\Lambda_{k},d_{j})$, $i,j<k$ are
isomorphic as well. Moreover, there exists an isomorphism of the complexes
$(\Lambda_{k},d_{k-1})$ and $(\Lambda_{k},d_{k})$. Therefore $(\Lambda
_{k},d_{i})$ and $(\Lambda_{k},d_{j})$ are isomorphic for any $i,j\leq k$. By
abusing the notation we again denote by $\kappa_{i,j}$ such isomorphism. Since
for any algebra automorphism which commutes with a given derivation $X$, the
corresponding automorphism of the differential form algebra commutes with the
Lie derivative along $X$, it is clear that $\kappa_{i,j}$ commutes with the
differentials $d_{l}$, $l\neq i,j$.\end{proof}

\begin{corollary}
\label{PropKappa1}For any permutation $\sigma\in S_{\mathbb{N}}$ there exists
an isomorphism $\kappa_{\sigma}$ of the multiple complexes $(\Lambda_{\infty
},d_{1},\ldots,d_{k},\ldots)$ and $(\Lambda_{\infty},d_{\sigma(1)}%
,\ldots,d_{\sigma(k)},\ldots)$.
\end{corollary}

\begin{proof}
Let $\sigma\in S_{\mathbb{N}}$ and $\Omega\in\Lambda_{k}\subset\Lambda
_{\infty}$. Let $l=\max\{\sigma(1),\ldots,\sigma(k)\}$ and $\tilde{\sigma}\in
S_{l}$ be such that $\tilde{\sigma}(i)=\sigma(i)$ for any $i\leq k$. Put
$\kappa_{\sigma}(\Omega)\equiv\kappa_{\tilde{\sigma}}(\Omega)$. $\kappa
_{\sigma}$ does not depend on the choice of $\tilde{\sigma}$ and is therefore
well-defined. Moreover, $\kappa_{\sigma}$ has inverse $\kappa_{\sigma^{-1}}$.
Thus it is an isomorphism of the multiple complexes $(\Lambda_{\infty}%
,d_{1},\ldots,d_{k},\ldots)$ and $(\Lambda_{\infty},d_{\sigma(1)}%
,\ldots,d_{\sigma(k)},\ldots)$.\end{proof}

The correspondence $\sigma\longmapsto\kappa_{\sigma}$ defines an action of the
group $S_{\mathbb{N}}$ on $\Lambda_{\infty}$.

An important fact is that the cohomology of complexes $(\Lambda_{\infty}%
,d_{k})$, $k\in\mathbb{N}$ is \textquotedblleft constant\textquotedblright.

\begin{theorem}
Complexes $(\Lambda_{2},d_{2})$ and $(\Lambda,d)$ are naturally homotopy equivalent.
\end{theorem}

\begin{proof}
In the proof of proposition \ref{PropKappa} a natural isomorphism
between the complexes $(\Lambda_{2}^{(0,\ast)},d_{2})$ and $(\Lambda,d)$ was
established. Let $\iota:(\Lambda_{2}^{(0,\ast)},d_{2})\longrightarrow
(\Lambda_{2},d_{2})$ be a natural embedding and $\pi:(\Lambda_{2}%
,d_{2})\longrightarrow(\Lambda_{2}^{(0,\ast)},d_{2})$ a natural projection. We
shall prove that the pair $(\iota,\pi)$ is a homotopy equivalence. Since
$\pi\circ\iota=\mathrm{id}_{\Lambda_{2}^{(0,\ast)}}$ it suffices to prove that
$\iota\circ\pi$ is homotopy equivalent to $\mathrm{id}_{\Lambda_{2}}$. If
$\Omega\in\Lambda_{2}^{(i,j)}$, then
\[
(\mathrm{id}_{\Lambda_{2}}-\iota\circ\pi)(\Omega)=\left\{
\begin{array}
[c]{cc}%
\Omega & \text{if }i\neq0\\
0 & \text{if }i=0
\end{array}
\right.  .
\]
Denote by $d^{0}$ the $0$--degree component of the ordinary de Rham
differential. It is a $\Lambda^{1}$--valued derivation of $A$. The
corresponding insertion operator $C=i_{d^{0}}: \Lambda\longrightarrow\Lambda$
is a derivation of $\Lambda$ and we have $C(\sigma)=s\sigma$ for $\sigma
\in\Lambda^{s}$.

Consider the $0$--degree map $i_{C}^{(2)}:\Lambda_{2}\longrightarrow
\Lambda_{2}$ and define the map $H_{2}:\Lambda_{2}\longrightarrow\Lambda_{2}$
by:
\[
H_{2}(\Omega)=\left\{
\begin{array}
[c]{cc}%
\tfrac{1}{s}i_{C}^{(2)}(\Omega) & \text{if }s\neq0\\
0 & \text{if }s=0
\end{array}
\right.  ,
\]
for $\Omega\in\Lambda_{2}^{(s,t)}$. Prove that $H_{2}$ is a homotopy
connecting the chain maps $\iota$ and $\pi$. Indeed, if $s=0$, then
$[H_{2},d_{2}](\Omega)=0=(\mathrm{id}_{\Lambda_{2}}-\iota\circ\pi)(\Omega).$
If $s\neq0$, then $[H_{2},d_{2}](\Omega)=\tfrac{1}{s}[i_{C}^{(2)}%
,d_{2}](\Omega)=\tfrac{1}{s}\mathcal{L}_{C}^{(2)}\Omega$. Show that $\mathcal{L}_{C}^{(2)}\Omega=s\Omega$. Assume that $\Omega=\sum
\sigma_{\alpha_{1} \ldots\alpha_{l}}\wedge d_{2}\sigma_{\alpha_{1}}%
\wedge\cdots\wedge d_{2} \sigma_{\alpha_{l}}$ for $\sigma_{\alpha_{1}%
\cdots\alpha_{l}},\sigma_{\alpha_{1}},\ldots,\sigma_{\alpha_{l}}\in\Lambda$.
Then
\begin{align*}
\mathcal{L}_{C}\Omega   =&\sum C(\sigma_{\alpha_{1}\cdots\alpha_{l}})\wedge
d_{2}\sigma_{\alpha_{1}}\wedge\cdots\wedge d_{2}\sigma_{\alpha_{l}}+\sum
\sigma_{\alpha_{1}\cdots\alpha_{l}}\wedge d_{2}C(\sigma_{\alpha_{1}}%
)\wedge\sigma_{\alpha_2}\wedge\cdots\wedge d_{2}\sigma_{\alpha_{l}}\\
  &+\cdots+\sum\sigma_{\alpha_{1}\cdots\alpha_{l}}\wedge d_{2}\sigma_{\alpha_{1}}\wedge\cdots\wedge d_{2}\sigma_{\alpha_{l-1}}\wedge d_{2}C(\sigma_{\alpha_{l}})=s\Omega.
\end{align*}
Therefore $[H_{2},d_{2}](\Omega)=\Omega=(\mathrm{id}_{\Lambda_{2}}-\iota
\circ\pi)(\Omega)$.\end{proof}

The following consequence is obvious.

\begin{corollary}
for any $k$ $H(\Lambda_{\infty},d_{k})\simeq H(\Lambda,d)$.
\end{corollary}

\section{Iterated differential forms over a smooth manifold}

In this and the next two sections our considerations are restricted to the
case $A=C^{\infty}(M)$, where $M$ is an $n$--dimensional smooth manifold. In
this situation it is convenient to pass to the category $\mathbf{Mod}_{A}^{g}$
of geometric $A$--modules (see, e.g., \cite{n03}) and to work with
representing objects of functors of differential calculus only in this
category. For instance, geometric differential forms over the algebra $A$ are
nothing but standard differential forms over the manifold $M$.

Denote by $\Lambda_{\infty}(M)$ (or simply $\Lambda_{\infty}$) iterated
geometric differential forms over $C^{\infty}(M)$. Their coordinate
description is as follows. Let $K=\{k_{1},\ldots,k_{r}\}\subset\mathbb{N}$ be
a finite subset. For $f\in A$ put $d_{K}f\overset{\mathrm{def}}{=} d_{k_{1}%
}\cdots d_{k_{r}}f$. If $x^{1},\ldots,x^{n}$ are local coordinates on $M$,
then $\Lambda_{\infty}(M)$ is locally generated as an algebra by the elements
$d_{K}x^{\mu}$, $\mu=1,\ldots,n$, $K\subset\mathbb{N}$. In particular, it is
easy to see that locally
\[
d_{K} f = \sum_{\{J_{1},\ldots,J_{l}\}}\frac{\partial^{l}f}{\partial
x^{\mu_{1}}\cdots\partial x^{\mu_{l}}}d_{J_{1}}x^{\mu_{1}}\wedge\cdots\wedge
d_{J_{l}}x^{\mu_{l}},
\]
where the sum runs over all repeated indexes and all partitions $\{J_{1}%
,\ldots,J_{l}\}$ of $K$ into $l$ parts, $1\leq l\leq r$.

Let $N$ be an $m$--dimensional manifold, $(y^{1},\ldots,y^{m})$ local
coordinates on $N$ and $\phi:M\longrightarrow N$, $x^{\mu}\longmapsto
y^{\alpha}=\phi^{\alpha}(x)$, $\mu=1,\ldots,n$, $\alpha=1,\ldots,m$, be a
smooth map. Denote by $\Lambda_{\infty}^{\prime}$ the algebra of iterated
differential forms over $N$ and consider the homomorphism $\Lambda_{\infty
}(\phi^{\ast}):\Lambda_{\infty}^{\prime}\longrightarrow\Lambda_{\infty}$
generated by $\phi^{\ast}:C^{\infty}(N)\longrightarrow C^{\infty}(M)$ which,
for simplicity, will be also denoted by $\phi^{\ast}$ as for ordinary forms.
Then, locally
\[
\phi^{\ast}(d_{K}y^{\alpha})=d_{K}(\phi^{\ast}(y^{\alpha}))=\sum
_{\{J_{1},\ldots,J_{l}\}}\frac{\partial^{l}\phi^{\alpha}}{\partial x^{\mu_{1}%
}\cdots\partial x^{\mu_{l}}}d_{J_{1}}x^{\mu_{1}}\wedge\cdots\wedge d_{J_{l}%
}x^{\mu_{l}}\in\Lambda_{\infty},
\]
where the sum runs over all repeated indexes and all partitions $\{J_{1}%
,\ldots,J_{l}\}$ of $K=\{k_{1},\ldots,k_{r}\}\subset\mathbb{N}$ into $l$
parts, $1\leq l\leq r$.

\section{Covariant tensors as iterated differential forms}

For any $p\in\mathbb{N}$ the functor $\mathrm{D}_{A}^{p}\equiv\mathrm{D}%
_{A}\circ\cdots\circ\mathrm{D}_{A}:\mathbf{Mod}_{A}^{g}\longrightarrow
\mathbf{Mod}_{A}^{g}$ is represented by the module $T_{p}^{0}(M)\equiv
\Lambda^{1}(M)^{\otimes p}$ of covariant $p$--tensors on $M$. For any $p$ the
map
\[
\square:A\times\cdots\times A\ni(f_{1},\ldots,f_{p})\longmapsto d_{1}%
f_{1}\wedge\cdots\wedge d_{p}f_{p}\in\Lambda_{\infty}%
\]
is a multi--derivation, i.e. $\square\in$ $\mathrm{D}_{A}^{p}(\Lambda_{\infty
})$. Therefore, there exists a unique $A$--homomorphism $\iota_{p}:T_{p}%
^{0}(M)\longrightarrow\Lambda_{\infty}$ such that
\[
\iota_{p}(df_{1}\otimes\cdots\otimes df_{p})=\square(f_{1},\ldots,f_{p}%
)=d_{1}f_{1}\wedge\cdots\wedge d_{p}f_{p}\in\Lambda_{\infty}%
\]
for any $f_{1},\ldots,f_{p}\in A$.

\begin{proposition}
\label{PropIota}$\iota_{p}$ is injective.
\end{proposition}

\begin{proof}
Obvious from local expression.\end{proof}

Proposition \ref{PropIota} shows that the calculus of covariant tensors over
$M$ is just a part of differential calculus over the algebra $\Lambda_{\infty
}$ and, therefore, is not \emph{conceptually closed}. In that sense the
proposed embedding of tensors into $\Lambda_{\infty}$ may be seen as a
conceptual closure of tensor calculus.

Standard operations with tensors, such as multiplications, \textquotedblleft
permutations of indexes\textquotedblright, insertions of vector fields, Lie
derivatives, etc, have proper counterparts in $\Lambda_{\infty}$ and this
looks as follows.

First, all kinds of tensor multiplications are encoded in the wedge product
\textquotedblleft$\wedge$\textquotedblright\ in $\Lambda_{\infty}$. For
instance, exterior and symmetric products of differentials $df$ and $dg$ look
as $d_{1}f\wedge d_{2}g-d_{1}g\wedge d_{2}f$ and $d_{1}f\wedge d_{2}%
g+d_{1}g\wedge d_{2}f$, respectively. Second, the embedding $\iota_{p}$ is
equivariant with respect to the natural action $\tau_{p}$ of the permutation
group $S_{p}$ on $T_{p}^{0}(M)$, i.e., $\iota_{p}\circ\tau_{p}(\sigma
)=\kappa_{\sigma}\circ\iota_{p}$, $\forall\sigma\in S_{p}$.

Evaluation of a tensor field $T\in T_{p}^{0}(M)$ on a $p$--ple of vector
fields $X_{1},\ldots,X_{p}\in\mathrm{D}(M)$ is interpreted in $\Lambda
_{\infty}$ by means of the formula
\begin{equation}
T(X_{1},\ldots,X_{p}) =(i_{X_{p}}^{(p)}\circ\cdots\circ i_{X_{1}}^{(1)}%
)(\iota_{p}(T))\in A\subset\Lambda_{\infty}. \label{Eq2}%
\end{equation}
Indeed, if $T=df_{1}\otimes\cdots\otimes df_{p}$, $f_{1},\ldots,f_{p}\in
C^{\infty}(M)$, then
\begin{align*}
(i_{X_{p}}^{(p)}\circ\cdots\circ i_{X_{1}}^{(1)})(\iota_{p}(T))  &
=(i_{X_{p}}^{(p)}\circ\cdots\circ i_{X_{1}}^{(1)})(d_{1}f_{1}\wedge
\cdots\wedge d_{p}f_{p}) =X_{1}(f_{1})\cdot\cdots\cdot X_{p}(f_{p})\\
&  =T(X_{1},\ldots,X_{p}).
\end{align*}
In particular, the insertion of a vector field $X\in\mathrm{D}(M)$ into the
$l$--th place of $T$ is given by
\[
\iota_{p}(T({}\cdot{},\ldots,{}\cdot{},\underset{l}{X},{}\cdot{},\ldots
,{}\cdot{}))=(i_{X}^{(l)}\circ\iota_{p})(T).
\]

Similarly, the Lie derivative $\mathcal{L}_{X}T$ of a tensor field $T\in
T_{p}^{0}(M)$ along a vector field $X\in\mathrm{D}(M)$ is given by
\[
\iota_{p}(\mathcal{L}_{X}T)=(X\circ\iota_{p})(T).
\]
Indeed, let $T=\sum g_{\alpha_{1}\cdots\alpha_{p}}df_{\alpha_{1}}\otimes
\cdots\otimes df_{\alpha_{p}}$, $f_{1},\ldots,f_{p}\in C^{\infty}(M)$. Then
\begin{align*}
(X\circ\iota_{p})(T)  &  =X({\textstyle\sum}g_{\alpha_{1}\cdots\alpha_{p}%
}d_{1}f_{\alpha_{1}}\wedge\cdots\wedge d_{p}f_{\alpha_{p}})={\textstyle\sum
}X(g_{\alpha_{1}\cdots\alpha_{p}})d_{1}f_{\alpha_{1}}\wedge\cdots\wedge
d_{p}f_{\alpha_{p}}\\
&  +{\textstyle\sum}g_{\alpha_{1}\cdots\alpha_{p}}d_{1}X(f_{\alpha_{1}}%
)\wedge\cdots\wedge d_{p}f_{\alpha_{p}}+\cdots+{\textstyle\sum}g_{\alpha
_{1}\cdots\alpha_{p}}d_{1}f_{\alpha_{1}}\wedge\cdots\wedge d_{p}%
X(f_{\alpha_{p}})\\
&  =\iota_{p}({\textstyle\sum}X(g_{\alpha_{1}\cdots\alpha_{p}})df_{\alpha_{1}%
}\otimes\cdots\otimes df_{\alpha_{p}}+{\textstyle\sum}g_{\alpha_{1}%
\cdots\alpha_{p}}dX(f_{\alpha_{1}})\otimes\cdots\otimes df_{\alpha_{p}}\\
&  +\cdots+{\textstyle\sum}g_{\alpha_{1}\cdots\alpha_{p}}df_{\alpha_{1}%
}\otimes\cdots\otimes dX(f_{\alpha_{p}}))\\
&  =\iota_{p}(\mathcal{L}_{X}T).
\end{align*}

The intrinsic characterization of covariant tensors as elements of
$\Lambda_{\infty}$ is as follows.

\begin{proposition}
\label{Tensors}A homogeneous element $\Omega\in\Lambda_{p}\subset
\Lambda_{\infty}$ of multi--degree $(1,1,\ldots,1)\in\mathbb{Z}^{p}$ is a
covariant tensor on $M$, i.e., $\Omega\in\mathrm{im}\ \iota_{p}$, iff the map
\[
\tilde{\Omega}:\mathrm{D}(M)\times\cdots\times\mathrm{D}(M)\ni(X_{1}%
,\ldots,X_{p})\longmapsto(i_{X_{p}}^{(p)}\circ\cdots\circ i_{X_{1}}%
^{(1)})(\Omega)\in A\subset\Lambda_{\infty}%
\]
is $A$--multi--linear.
\end{proposition}

\begin{proof}
Let $\Omega\in\Lambda_{p}$ be a homogeneous element of multi--degree
$(1,1,\ldots,1)$. $\Omega$ can be expressed in the form
\[
\Omega=\iota_{p}(T)+\sum_{\{J_{1},\ldots,J_{l}\}}g_{\alpha_{1}\cdots\alpha
_{p}}^{J_{1}\cdots J_{l}}d_{J_{1}}(f_{\alpha_{1}J_{1}})\wedge\cdots\wedge
d_{J_{l}}(f_{\alpha_{l}J_{l}})
\]
where $T\in T_{p}^{0}(M)$, $g_{\alpha_{1}\cdots\alpha_{p}}^{J_{1}\cdots J_{l}%
},f_{\alpha_{1}J_{1}},\ldots,f_{\alpha_{l}J_{l}}\in C^{\infty}(M)$ and the sum
runs over $\alpha_{1},\ldots,\alpha_{l}$ and all partitions $\{J_{1}%
,\ldots,J_{l}\}$ of $\{1,\ldots,p\}$ into $l$ parts, $1\leq l<p$. Let
$X\in\mathrm{D}(M)$ and $s\leq p$. Without loss of generality suppose that
$J_{1}=\{i_{1},\ldots,i_{r},s\}$. Then,
\begin{align*}
i_{X}^{(s)}\Omega=  &  \iota_{p}(T({}\cdot{},\ldots,{}\cdot{},\underset{s}%
{X},{}\cdot{},\ldots,{}\cdot{}))\\
&  +\sum_{\{J_{1},\ldots,J_{l}\}}g_{\alpha_{1}\cdots\alpha_{p}}^{J_{1}\cdots
J_{l}}d_{\{i_{1},\ldots,i_{r}\}}(X(f_{\alpha_{1}J_{1}}))\wedge d_{J_{2}%
}(f_{\alpha_{2}J_{2}})\wedge\cdots\wedge d_{J_{l}}(f_{\alpha_{l}J_{l}}).
\end{align*}
Thus, $\tilde{\Omega}$ is multilinear iff $\sum_{\{J_{1},\ldots,J_{l}%
\}}g_{\alpha_{1}\cdots\alpha_{p}}^{J_{1}\cdots J_{l}}d_{J_{1}}(f_{\alpha
_{1}J_{1}})\wedge\cdots\wedge d_{J_{l}}(f_{\alpha_{l}J_{l}})=0$, i.e.
$\Omega=\iota_{p}(T)$.
\end{proof}


\begin{thebibliography}{9}                                                                                                %


\bibitem{v98}A.~M.~Vinogradov, in M.~Henneaux, I.~Krasil'shchik, and
A.~Vinogradov (Eds.), \emph{Secondary Calculus and Cohomological Physics},
Contemporary Mathematics \textbf{219}, AMS (1998). See also The Diffiety
Inst.~Preprint Series, DIPS 5/98.

\bibitem{v01}A.~M.~Vinogradov, \emph{Cohomological Analysis of Partial
Differential Equations and Secondary Calculus}, AMS \textquotedblleft
Translation of Mathematical Monographs \textquotedblright series \textbf{204} (2001).

\bibitem{v72}A.~M.~Vinogradov, \emph{Soviet Math.~Dokl.~}\textbf{13} (1972) 1058.

\bibitem{v81}A.~M.~Vinogradov, \emph{J.~Soviet Math.~}\textbf{17} (1981) 1624.

\bibitem{v89}M.~M.~Vinogradov, \emph{Russian Math.~Surveys }\textbf{44}, n${}^\circ$ 3
(1989) 220.

\bibitem{v96}A.~Verbovetsky, \emph{J.~Geom.~Phys.~}\textbf{18} (1996) 195.

\bibitem{vv06}A.~M.~Vinogradov and L.~Vitagliano, \emph{Dokl.~Math.~}\textbf{73}, n${}^\circ$ 2 (2006) 169.

\bibitem{n03}J.~Nestruev, \emph{Smooth Manifolds and Observables},
Springer--Verlag (New York) 2003.

\end{thebibliography}
\end{document}